\documentclass[12pt]{amsart}
\usepackage[colorlinks=true, pdfstartview=FitV, linkcolor=blue, citecolor=blue, urlcolor=blue]{hyperref}

\usepackage{amssymb,amsmath, amscd}
\usepackage{times, verbatim}
\usepackage{graphicx}
\usepackage[english]{babel}
 \usepackage[usenames, dvipsnames]{color}
\usepackage{amsmath,amssymb,amsfonts}
\usepackage{enumerate}

\usepackage{tikz}
\usetikzlibrary{calc, arrows, decorations.markings} \tikzset{>=latex}
\usepackage{amsmath,amssymb,amsfonts}
\usepackage{amscd, amssymb, latexsym, amsmath, amscd}
\usepackage[all]{xy}
\usepackage{pb-diagram}
\usepackage{enumerate}
\usepackage{ulem}
\usepackage{anysize}
\marginsize{3cm}{3cm}{3cm}{3cm}
\input xy
\xyoption{all}
\usepackage{pb-diagram}
\usepackage[all]{xy}
\input xy
\xyoption{all}

\theoremstyle{plain}

\theoremstyle{definition}
\newtheorem{Question}{\bf Question}
\newtheorem{Remark}{Remark}

\theoremstyle{definition} \theoremstyle{definition}

\theoremstyle{remark}

\newcommand{\G}{\textsc{\G}}

\newcommand{\C}{\mathbb{C}}
\newcommand{\norm}{\mathbb{N}}
\newcommand{\Fq}{\mathbb{F}_q}
\newcommand{\F}{\mathbb{F}}

\def\G{{\rm G}}

\def\SL{{\rm SL}}

\def\GL{{\rm GL}}

\def\OF{{\rm O}_F}

\subjclass{Primary 11F70; Secondary 22E55}

\begin{document}

\title{Some questions about representations of p-adic groups}
\author{Dipendra Prasad}
\thanks{I thank Prof Kumar Balasubramanian, and IISER Bhopal, for hosting a discussion meeting on the subject in July 2025, and for the excellent hospitality, cf. here:  \href{https://sites.google.com/iiserb.ac.in/discussion-meeting-2025/home}.
 Thanks are also due to Anand Chitrao for making a first draft of these notes, and to various participants for their interest.}

\begin{abstract} Some question about representations of $p$-adic groups
  are discussed.
\end{abstract}

  \address{Indian Institute of Technology Bombay, Powai, Mumbai-400076}
\email{prasad.dipendra@gmail.com}
\date{\today}

\maketitle

\vspace{1cm}

\tableofcontents

\vspace{1cm}

During a
discussion meeting on representation theory held at IISER,
Bhopal in July 2025,
there was a request from several participants if I could suggest
some unsolved but approachable questions without too much background. I gave the last two lectures of the program discussing some of these questions which appear here,
together with some  background material with each question. The questions
posed here are on themes related to the topics which came up in the discussion meeting, and thus are to be viewed from a very narrow lens. 

Throughout this document, we let $F$ denote a local field with the ring of integers $\OF$, uniformizer $\varpi$ and residue field $\Fq$.
The non-archimedean absolute value of an element $x \in F$ will be denoted by
$\norm{x}$, and $\nu$ the character of $F^\times$ given by $\nu(x)=\norm(x)$.

\section{Representations of $p$-adic group containing a representation of the group over the finite field}
Let $\G$ be a split group, or more generally an unramified group over $F$, thus a quasi-split group over $F$
which splits over an unramified extension of $F$. The group $\G(\OF)$ is then a maximal compact subgroup in $\G(F)$ and comes equipped with a
surjective map onto the finite group $\G(\Fq)$. Therefore, one may think of representations of $\G(\Fq)$ as representations of $\G(\OF)$ by inflation.
In particular, inflating a cuspidal representation  $\rho$ of $\G(\Fq)$ to $\G(\OF)$ and then inducing to $\G(F)$, one gets an irreducible
supercuspidal representation $\mathrm{ind}_{\G(\OF)}^{\G(F)}\rho$ of $\G(F)$ of depth $0$ if the connected center of $\G$ is trivial.

    \vspace{0.5cm}
    
    \begin{Question}
      Classify all
    irreducible representations of $\G(F)$ containing the Steinberg representation of $\G(\Fq)$ as a $\G(\OF)$-subrepresentation.
    \end{Question}
    
    \vspace{0.3cm}

    \noindent \textbf{Expectation:}
    Irreducible representations of $\G(F)$ containing the Steinberg representation of $\G(\Fq)$ as a $\G(\OF)$-subrepresentation are:

    \begin{enumerate}
    \item a twist of the Steinberg representation of $\G(F)$ by a character of $\G(F)$, or,
    \item  generic, unramified representations of $\G(F)$, thus by a well-known result due to Jian-Shu Li, \cite{Li}, these are irreducible spherical
      principal series representations of $\G(F)$.
    \end{enumerate}
    
    \vspace {0.5cm}

Here is a more general question.    Let $\pi_\lambda$ be the irreducible representation of $\G(\Fq)$ associated to an irreducible  representation $\lambda$ of the Weyl group of $\G$, appearing inside ${\rm Ind}_{B(\Fq)}^{\G(\Fq)} (\C)$.

    \vspace{0.5cm}
    
    \begin{Question}

    Classify representations of $\G(F)$ containing $\pi_\lambda$. This classification can be in terms of the (enhanced) Langlands parameter, or for the case of $\mathrm{GL}_n$, in terms of the Zelevinsky classification.
    \end{Question}
    
    \vspace{0.5cm}

    Associated with each standard parabolic $P$, there is an irreducible representation $\Pi_P$ on the space of locally constant functions
    on $\G/P$ modulo functions on $\G/Q$ for
    parabolics $Q$ strictly containing $P$. (It is a theorem of Casselman that $\pi_P$ is an irreducible representation of $\G(F)$ when $F$ is a p-adic field. However, when $F=\Fq$ is a finite field, $\pi_P$ is typically not irreducible.)
A very special case of the question is to classify $(P,\lambda)$ for which
$\mathrm{Hom}_{\mathrm{GL}_n(\OF)}(\Pi_P, \pi_{\lambda}) \neq 0$? For the case
of $\GL_n(F)$,
    one may guess that
    $\mathrm{Hom}_{\mathrm{GL}_n(\OF)}(\Pi_P, \pi_{\lambda}) \neq 0$
    if and only if  $P \mapsto \lambda$, under the natural surjective map from the set of standard parabolics $P$ of $\mathrm{GL}_n(F)$ to the set of partitions $\lambda$ of $n$.

    \vspace{1cm}

    \begin{Remark}
      Representations considered in this section have minimal ramification, in particular, they have fixed vectors under
    an Iwahori subgroup, thus by a famous theorem of Borel, cf. \cite{Bo}, they arise as subquotients of unramified principal series.
But even here, it appears to me that the
    precise relationship as formulated in the questions above,  between enhanced Langlands parameter for $\G(F)$
    with Lusztig's works classifying representations of $\G(\Fq)$, is not clarified.
    \end{Remark}
\newpage
    
\section{Degenerate Whittaker models of representations of $\GL_n(D)$}
Let $D$ be a division algebra. Fix a non-trivial additive
character $\phi : D \to \mathbb{C}^\times$. A representation $\Pi$ of $\mathrm{GL}_n(D)$ is said to have
a degenerate Whittaker model if the twisted Jacquet module $\Pi_{N, \psi} \neq 0$, where $N$ is the subgroup of upper triangular unipotent matrices in $\mathrm{GL}_n(D)$ and
    \[
        \psi \begin{pmatrix} 1 & X_1 & * & \cdots & * \\ 0 & 1 & X_2 & \ddots & * \\ \vdots & \vdots & \ddots & \ddots & \vdots \\ 0 & 0 & \cdots & 1 & X_{n - 1} \\ 0 & 0 & \cdots & 0 & 1 \end{pmatrix} = \phi(\mathrm{tr}(X_1 + X_2 + \cdots + X_{n - 1})).
    \]

    The twisted  Jacquet module $\Pi_{N, \psi} $ is a representation
    of $D^\times$, and its dimension is an important invariant of the representation $\pi$, which by the work of Moeglin-Waldspurger, \cite{MW}, is known to be finite (though, strictly speaking, perhaps only in characteristic zero).
    If $D=F$,   the twisted  Jacquet module $\Pi_{N, \psi} $ is just the Whittaker model of $\pi$.

    \begin{Question} Classify irreducible
      representations of $\mathrm{GL}_n(D)$ which have a non-zero degenerate Whittaker model. One could ask for a classification either in terms of
      the Jacquet-Langlands correspondence from
      irreducible representations of
      $\mathrm{GL}_n(D)$ to irreducible representations of
      $\GL_{mn}(F)$, or through
      the Langlands parameters.
            Note that the ``usual'' argument with
    the mirabolic in $\GL_n(F)$ for supercuspidals works here too to show that
    supercuspidals always have a degenerate Whittaker model for $\GL_n(D)$,
   and globally, cuspidal representations  of $\GL_n(D)$ have a global
   degenerate Whittaker model. Further, the Bernstein-Zelevinski theory of ``derivatives'', and the Leibnitz
   rule work just as well (except that the highest derivative  $\pi^n$  for $\GL_n(D)$ is not necessarily one dimensional). Therefore all tempered representations of $\GL_n(D)$
    have degenerate Whittaker model.
  Is there an analogue of the ``standard module conjecture''?
    Is there a way to determine the dimension
    of the space of degenerate Whittaker models
        in terms of adjoint $\gamma$-functions as in the Hiraga-Ichino-Ikeda
    conjecture in \cite{HII}?
      
\end{Question}

    \begin{Question} Classify irreducible
      representations of $\mathrm{GL}_n(D)$ with a unique degenerate Whittaker model. In particular, which cuspidal representations have this property?
\end{Question}
    \vspace{0.2cm}
    
    \begin{Remark}
    The above question is there because I do not know how to answer
    it even in the simplest case of $\GL_2(D)$, $D$ quaternion.   Presumably, the only cuspidal representations
 of $\GL_n(D)$, $D$ any division algebra, 
    with one dimensional
    space of degenerate Whittaker models are those of depth zero?
    (That the cuspidal representations of depth zero have a
   unique degenerate Whittaker model follows from the uniqueness of
    Whittaker models over finite field.)
    More generally, one would like to understand
    how the dimension of the space of degenerate Whittaker models  grow
    with the conductor of the representation?
    \end{Remark}
    
    \begin{Remark} Questions on degenerate Whittaker model
        formulated above for $\GL_n(D)$
 could also be formulated for general
    $p$-adic reductive group, taking $P=MAN$ to be a minimal parabolic, and
    taking a character $\psi:N \rightarrow \C^\times $ which is non-trivial on all simple roots
    of $A$ on $N$. The question then is to classify representations
    $\pi$ of $\G$ with $\pi_{N,\psi} \not = 0$.
    \end{Remark}
    \newpage
    
    \section{Twisted-Jacquet modules for $\mathrm{GL}_{2n}(\Fq)$} 
    Let $\pi$ be a cuspidal representation of $\mathrm{GL}_n(\Fq)$. It is well-known that the principal series $\pi \times \pi$ sits in the following exact sequence
    \[0 \to \mathrm{St}_2(\pi) \to \pi \times \pi \to \mathrm{Sp}_2(\pi) \to 0,        \] where  $\mathrm{St}_2(\pi)$ is the generic component
        of $\pi \times \pi$, that we call generalised Steinberg representation
        of $\GL_{2n}(\Fq)$, and we call  $\mathrm{Sp}_2(\pi)$ as the generalized
        trivial representation   of $\GL_{2n}(\Fq)$. One knows that,
        \[ \dim \mathrm{St}_2(\pi) = q^n \dim \mathrm{Sp}_2(\pi).\]

        We will be considering the twisted Jacquet modules for
        the character $\psi$ of the unipotent radical $N = M(n,\Fq)$ of the $(n, n)$-parabolic given by $x\rightarrow \phi(tr(X))$ where $\phi$ is any fixed non-trivial character of $\Fq$.
        A calculation of this was done for supercuspidal representations of
        $\GL_{2n}(\Fq)$ in \cite{Pr}. It was also mentioned there that the
        same calculation works for the Deligne-Lusztig induction $R(T,\theta)$
        where $T = \F_{q^{2n}}^\times$, and $\theta: T = \F_{q^{2n}}^\times \rightarrow \mathbb{C}^\times$ to be any character, and not necessarily one which
        gives rise to a cuspidal representation of $\GL_{2n}(\Fq)$.

        Suppose $\pi = R(\F_{q^{n}}^\times, \theta_0)$, a cuspidal representation of $\GL_n(\Fq)$, where
        $\theta_0: \F_{q^{n}}^\times \rightarrow \mathbb{C}^\times$ is a
        non-degenerate character. The character $\theta_0$ on composition with the norm mapping from $\F_{q^{2n}}^\times$ gives rise to a character, call it $\theta$, of  $\F_{q^{2n}}^\times$, and hence the Deligne-Lusztig representation
        $R(T,\theta)$, for $T = \F_{q^{2n}}^\times$.
        One knows that,
        \[ R(T,\theta) = \mathrm{St}_2(\pi)- \mathrm{Sp}_2(\pi).\]

        As the results of \cite{Pr} about twisted Jacquet module
        are valid for $R(T,\theta)_{N, \psi}$, for any character  $\theta: T =
        \F_{q^{2n}}^\times \rightarrow \mathbb{C}^\times$, and not necessarily one which
        gives rise to a cuspidal representation of $\GL_{2n}(\Fq)$, we find that:

    \[
    \mathrm{Ind}_{\mathbb{F}_{q^{n}}^*}^{\mathrm{GL}_n(\Fq)} \theta^2 =
    \mathrm{St}_2(\pi)_{N, \psi} - \mathrm{Sp}_2(\pi)_{N, \psi}.
    \]
    Also,
    \[
    (\pi \times \pi)_{N, \psi} = \mathrm{St}_2(\pi)_{N, \psi} + \mathrm{Sp}_2(\pi)_{N, \psi}.
       \]

Combining these two equations, one obtains 
\[\mathrm{Sp}_2(\pi)_{N, \psi} = \frac{1}{2}\left((\pi \times \pi)_{N, \psi}
- \mathrm{Ind}_{\mathbb{F}_{q^n}^*}^{\mathrm{GL}_n(\Fq)}\theta^2\right)\]

Now, as $(\pi \times \pi)_{N, \psi}$ can be calculated via a direct application of Mackey theory, and was done in \cite{Pr} for $\pi$ a cuspidal representation of $\GL_2(\F_q)$, one can hope to find an explicit answer to the structure of $\mathrm{Sp}_2(\pi)_{N, \psi}$ as a module for $\GL_n(\F_q)$.
For the case of $\pi$,    a cuspidal representation of $\GL_2(\F_q)$, it is an easy consequence
of results proved in \cite{Pr} that,
$ \mathrm{Sp}_2(\pi)_{N, \psi} $ is the one dimensional representation $\omega_\pi$ of $\GL_2(\F_q)$. We pose the general question as follows:
    \vspace{0.5cm}
    
    \begin{Question} Let $\pi$ be a cuspidal representation of $\mathrm{GL}_n(\Fq)$. For the principal series representation $\pi \times \pi$ of $\GL_{2n}(\F_q)$, calculate
      its twisted Jacquet module $(\pi \times \pi)_{N, \psi}$, and the twisted Jacquet module $\mathrm{Sp}_2(\pi)_{N, \psi}$ of the corresponding ``generalised trivial'' representation. 
    \end{Question}

    \vspace{.5 cm}

    \begin{Remark} The question above will allow one to come up with an answer
    to a similar question on $\GL_2(D)$, where $D$ is any division
    algebra over $p$-adics of dimension $n^2$, for the principal series $\nu^{a} \pi \times \nu^{-a} \pi
    $ of $\GL_2(D)$ which sits in an exact sequence:
    \[0 \to \mathrm{St}_2(\pi) \to \nu^{a} \pi \times \nu^{-a}\pi \to \mathrm{Sp}_2(\pi) \to 0,        \] where  $\pi$ is a finite dimensional irreducible
    representation of $D^\times$, and $2a$ is an integer $> 0$, depending on $\pi$, with $a = \frac{1}{2}$, if the JL transfer of $\pi$ to $\GL_n(F)$ is cuspidal. It is easy to see that  $(\nu^{a} \pi \times \nu^{-a}\pi)_{N,\psi} = \pi \otimes \pi$ as a representation of $D^\times$, but how it distributes itself in the two components,   $\mathrm{St}_2(\pi)_{N,\psi}$ and $  \mathrm{Sp}_2(\pi)_{N,\psi}$ is not understood at the moment, except when $D$ is a quaternion division algebra, cf. \cite{NS}.
    \end{Remark}

    \begin{Question}
      For $\pi = R(T,\theta)$, an irreducible cuspidal representation of
      $\GL_{2n}(\F_q)$, it was proved in \cite{Pr} that, in the notation
      above,
      \[\pi_{N,\psi} = {\rm Ind}_{\F_{q^n}^\times}^{\GL_n(\F_q}(\theta|_{\F_{q^n}^\times})
      \hspace{4cm}{(1)}.\]

      What's the analogue of this theorem for $p$-adic $\GL_{2n}(F)$?
      As ``most'' supercuspidal representations of $\GL_{2n}(F)$, $F$ a p-adic field, are ``induced'' from a character of a field extension of $F$ of degree $2n$, it is meaningful to ask if the assertion analogous to
      (1) above holds true for $F$ a p-adic field?
      However, as this result for finite fields was proved by explicit character theory, and as the number of irreducible components in $\pi_{N,\psi}$ is a polynomial in $q$ of degree $>0$, there cannot be a
      character theoretic proof for $F$ a $p$-adic field.
      One possibility for the formulation would be to assert that for $F$ a p-adic field, the two sides of the equation (1) have  finite filtrations  whose successive quotients are the same or are related. For this, one must find a more ``geometric'' proof of (1), and not character theoretic,
      in the case of finite fields! Is $\pi_{N,\psi}$, $F$ a p-adic field, of
      ``geometric'' origin, for $\pi$ supercuspidal in the same way as $\pi_{N,\psi}$ is of geometric origin when $\pi$ is a principal series representation?
      (As another example of a representation of infinite length of ``geometric origin'', is the Bernstein-Zelevinski filtration on a representation of $\GL_n(F)$ restricted to the mirabolic.) 
      \end{Question}

        \newpage
    
    \section{Covering groups}
    Let $\widetilde{\G} \rightarrow \G$
    be a finite cover of $p$-adic groups,
    a  central extension with $Z$ the kernel of the map 
    $\widetilde{\G} \rightarrow \G$ which we assume is a cyclic group, and has a
    fixed faithful character $\epsilon: Z\rightarrow \C^\times$.
    Let $K$ be a compact-modulo-center subgroup of $\G$ on which this cover splits. Fix a splitting $s : K \to \widetilde{\G}$. Let $\widetilde{K}= Z\cdot s(K)$.
    
    \vspace{0.5cm}
    
    \begin{Question}
      Suppose $\Pi = \mathrm{Ind}_{K}^{\G}(\rho)$ is an irreducible supercuspidal representation of $\G$ for $\rho$ an irreducible representation of $K$. Let
      $\widetilde{\rho}$ be the representation of  $\widetilde{K}= Z\cdot s(K)$ which
      is $\rho$ on $s(K)$, and $\epsilon$ on $Z$. Define
      $\widetilde{\Pi} = \mathrm{Ind}_{\widetilde{K}}^{\widetilde{\G}}(\widetilde{\rho})$. Show that
    \begin{enumerate}
        \item $\widetilde{\Pi}$ is irreducible and supercuspidal.
                \item Assume that  $\widetilde{\G} \rightarrow \G$ splits on all
          maximal compact-modulo-center subgroups of $\G$, as it does for $\G=\SL_n(F)$, $(p,d)=1$ where $p$ is the residue characteristic of $F$, and $d$ is the degree of the cover $\widetilde{\G} \rightarrow \G$. Is every supercuspidal
          of $\widetilde{\G}$ with central character $\epsilon$ on $Z$
          of the form as in (i) if every supercuspidal of $\G$ is induced from a compact-modulo-center subgroup of $\G$? If so, this will prove that all supercuspidal
of $\widetilde{\G}$ 
          are induced from a compact-modulo-center subgroup knowing it on $\G$.

\item $\widetilde{\Pi}$ is generic if and only if $\Pi$ is generic.
\item $\dim \mathrm{Hom}_N(\widetilde{\Pi}, \psi) = \dim \mathrm{Hom}_N({\Pi}, \psi)$, where $(N, \psi)$ is the Whittaker datum of $\G$. In particular,
  multipliplicity one holds for Whittaker models for these representations of the covering groups,
  unlike the higher multiplicity of \cite{GGK} for the principal series
  representations of covering groups.

    \end{enumerate}
\end{Question}
    \vspace{0.5cm}

    \begin{Remark} See the paper \cite{PP} for some related material. It would be nice to classify groups $\G$, together with the covering
      $\widetilde{\G} \rightarrow \G$ which split on all maximal compact-modulo-center subgroups of $\G$, which presumably happens in most cases.
      \end{Remark}
    \newpage
    
\section{Fourier transforms on $\mathrm{M}_n(\Fq)$}
It is well-known that $\mathrm{M}_n(F)$ has a rich
theory of the Fourier transform which plays a fundamental role in representation theory of $\GL_n(F)$ for which there is the very pleasant and seminal paper of Howe, cf. \cite{Ho}.
I have not seen any analogous works for finite fields, so here is a question for my own education which I
suspect is there in some form or the other in the existing literature.

    \vspace{0.5cm}
    
    \begin{Question}
Given a function $f$ on $\mathrm{M}_n(\F_q)$, its Fourier transform $\hat{f}$ is a function also on $\mathrm{M}_n(\Fq)$.
      Compute the Fourier transform of
    \begin{enumerate}
        \item the characteristic function of a nilpotent orbit (a conjugacy class containing nilpotent elements),
        \item the characteristic function of a semisimple orbit (a conjugacy class containing semisimple elements),
        \item the characteristic function of a nilpotent cone (union of all nilpotent orbits).
        \item Relate the above to character theory for $\GL_n(\F_q)$.
          \item Do the same for more general semi-simple Lie algebras $\mathfrak{g}$ over $\F_q$.
    \end{enumerate}
\end{Question}

\end{document}